\documentclass[12pt]{article}
\usepackage{latexsym, amsfonts, amsmath, amssymb, color, graphicx}
\usepackage{graphicx}

\newcommand{\qed}{\hskip 5mm \rule{2.5mm}{2.5mm}}
\newcommand{\R}{{\mathbb R}}

\newcommand{\N}{{\mathbb N}}

\newcommand{\proof}{{\em Proof:\ }}

\begin{document}
\newtheorem{thm}{Theorem}[section]
\newtheorem{defs}[thm]{Definition}
\newtheorem{lem}[thm]{Lemma}
\newtheorem{note}[thm]{Note}
\newtheorem{cor}[thm]{Corollary}
\newtheorem{prop}[thm]{Proposition}
\renewcommand{\theequation}{\arabic{section}.\arabic{equation}}
\newcommand{\newsection}[1]{\setcounter{equation}{0} \section{#1}}
\renewcommand{\baselinestretch}{1}
\title{Inverse problems with a general transfer condition 
      }
\author{Sonja Currie  \footnote{ Supported by NRF grant no. IFR2011040100017}\, $^{+}$\\
Marlena Nowaczyk \footnote{ Partially supported by Foundation for Polish Science, Programme Homing 2009/9} \,$^{o}$ \\
Bruce A. Watson \footnote{Supported by NRF grant no. IFR2011032400120} \, $^{+}$ \\ \\
\, $^{+}$ School of Mathematics\\
University of the Witwatersrand\\
Private Bag 3, P O WITS 2050, South Africa \\ \\
\, $^{o}$ AGH University of Science and Technology\\
Faculty of Applied Mathematics\\
al. A. Mickiewicza 30, 30-059 Krakow, Poland}
\maketitle
\abstract{
\noindent
We consider a Sturm-Liouville operator on a finite interval as well as a scattering problem on the real line both with transfer conditions at the origin. On a finite interval we show that the  the Titchmarsh-Weyl $m$-function can be uniquely determined from two spectra for the same equation but with varied boundary conditions at one end of the interval. In addition, we prove that the $m$-function can also be uniquely reconstructed from one spectrum and the corresponding norming constants. For the scattering problem on the real line we assume that the potential has compact essential support. For a given symmetric finite intervals containing the essential-support of the potential and a pair of separated boundary conditions imposed at the ends of the interval, the spectrum and corresponding norming constants can be uniquely recoverable from the scattering data on $\R$. Consequently the potential and transfer matrix can be determined.}

\parindent=0in
\parskip=.15in
\newsection{Introduction\label{sec-intro}}

Inverse spectral problems for the Sturm-Lioville operator were studied as early as 1929 by Ambartsumyan. The applications to mathematical physics lead to a renewed interest in the problem in recent years. For a comprehensive discussion on the early development of the inverse Sturm-Liouville problem see the introduction of Levitan's book \cite{Levitan}. More recently there has been interest in inverse Sturm-Liouville problems where the eigenparameter appears in the boundary conditions, for example see \cite{BBW} and the references therein. 

For an introduction to inverse scattering and inverse spectral problems including many direct applications we refer the reader to the book by Chadan et al. \cite{CCPR}. 
Inverse Sturm-Liouville problems with a discontinuity at an interior point were first 
studied by Hochstadt and Lieberman in \cite{hoch}. These results were then generalized in the 
famous paper by Hald, \cite{hald}, where it is shown that if the potential is known over half of the 
interval and if one of the boundary conditions is given, then the potential and the other boundary 
condition are uniquely determined by the eigenvalues. This in turn was extended by Willis in \cite{willis} to 
the case of two interior discontinuites. Similar techniques were then used by Kobayashi in \cite{koba} to give a 
uniqueness proof for the inverse Sturm-Liouville problem on a bounded interval with a symmetric 
potential having two interior jump discontinuities. A case with finitely many transmission conditions was discussed by 
Shahriari, Jodayree Akbarfam and Teschl in \cite{SJT}. In particular, the authors prove that the specification of the Weyl function and the weight function uniquely determines the Sturm-Liouville operator on a finite interval for both Robin and separated eigenparameter dependent boundary conditions. 

Ramm, \cite{ramm}, discusses inverse scattering and spectral problems 
on the half-line in detail.
 Some of the main topics included are, invertibility of the steps in the Gel'fand-Levitan and 
Marchenko inversion procedures, Krein inverse scattering theory and inverse problems. Aktosun and Weder, in \cite{Aktosun}, consider the Schr\"odinger equation on the half line with a real-valued,
integrable potential having a finite first moment. They extend the well-known two-spectrum uniqueness theorem of Borg
and Marchenko to the case where there is also a continuous spectrum. An elegant proof of the local Borg-Marchenko theorem is given by Bennewitz in \cite{Benne}.

In the last twenty years inverse problems on quantum graphs have received much attention. In particular Gerasimenko, in \cite{Gera}, develops a procedure for recovering the potential in the Schr\"odinger equation on a non-compact graph from the scattering data. Very recently, in \cite{Buterin}, using 
 the method of spectral mappings Buterin and Freiling prove that the specification of the spectral-scattering data uniquely determines the Sturm-Liouville operator on a non-compact graph. In addition, Dehghani and Jodayree Akbarfam study an inverse spectral problem on a three-star graph and show that if four particular spectra do not intersect then it is possible to recover the potential uniquely, see \cite{Deh}.
It should be noted here that the problem we consider in Section 3 is a special case of a non-compact graph consisting of two infinite edges connected at one vertex. 

Inverse scattering on the real line was discussed by Deift and Trubowitz in \cite{Deift} where they reconstruct a potential from its reflection coefficient, bound states and norming constants.
In \cite{alpay}, Alpay and Gohberg consider the inverse problem for Sturm-Liouville operators on the real line with rational reflection coefficient. They provide exact formulae for the potential associated with a Sturm-Liouville equation having rational reflection coefficient function.

In this paper we investigate the differential equation
\begin{equation}\label{diff}
 \ell y:=-\frac{d^2y}{dx^2} + q(x)y = {\zeta}^2 y,\quad \mbox{on } [-S,S],
\end{equation}
in $L^2[-S,S]$ with point transfer condition
\begin{equation}\label{tx-condition}
  \left[\begin{array}{c}y(0^+)\\ y'(0^+)\end{array}\right]=M
   \left[\begin{array}{c}y(0^-)\\ y'(0^-)\end{array}\right].
\end{equation}
Here the entries of $M$ are taken to be real, $q \in L^2[-S,S]$ is assumed to be real valued. We assume that $\det M >0$ and without loss of generality we will consider $\det M =1$. At the endpoints we impose the following boundary conditions
\begin{eqnarray}
y(-S)\cos \alpha - y'(-S)\sin \alpha &=0 \label{alpha},\quad\alpha\in [0,\pi),\\
y(S) \cos \beta -y'(S)\sin \beta &= 0 \label{beta},\quad\beta\in (0,\pi].
\end{eqnarray}

We will only consider point transfer matrices at the origin and henceforth will refer to 
them as transfer matrices. In the physical context the transfer matrix represents a change of 
medium which affects the incident wave as represented by 
components of the matrix. Our transfer matrices will be real 
constant transfer matrices i.e. all components will be constants. 

We show that two distinct problems are spectrally related. Firstly there is the finite interval problem (\ref{diff}) and (\ref{tx-condition}). For this problem we summarize, from \cite{cenw2}, how the Titchmarsh-Weyl $m$-function, $m(\lambda )$, of (\ref{diff}) on $[-S,S]$ obeying the transfer condition (\ref{tx-condition}) with separated boundary conditions at the end points is defined. An asymptotic approximation for $m(\lambda )$ is then obtained and consequently we show that $m(\lambda )$ can be uniquely determined from two spectra. In addition, we show that from one spectrum and norming constants it is also possible to uniquely find $m(\lambda )$. Secondly there is the scattering problem (\ref{diff}), (\ref{tx-condition}) on $(-\infty,0)\cup (0, \infty) $ where the potential $q$ has compact essential support in $[-S,S]$. In \cite{cenw2} we show that the scattering data uniquely determines the $m$-function, $m(\lambda )$, and indeed the converse holds i.e. given $m(\lambda )$ we can uniquely find the scattering data. In this paper we show that the scattering data also determines two spectra or one spectrum and norming constants. Thus the finite interval problem on $[-S,S]$ determines the scattering problem on the line with potential having compact essential support in $[-S,S]$ and vice-versa. 

The paper is organised as follows. Section $2$ is dedicated to asymptotic approximation for the $m$-function. Using these asymptotic estimates in Section 3 we prove that from two spectra the $m$-function can be uniquely determined. Moreover using residue calculus and the Mittag-Leffler expansion it is shown that $m$-function can also be reconstructed from one spectrum and the corresponding norming constants. Inverse scattering problems on the line with potential having compact essential support form the topic of Section $4.$ We prove that given the scattering data on $\R$, i.e. the reflection coefficient and eigenvalues, two spectra or one spectrum and corresponding norming constants on the finite interval $[-S,S]$ with $\makebox{ess sup } q \subset [-S,S]$ can be reconstructed. Consequently the operator and transfer matrix can be uniquely found.

\newsection{Titchmarsh-Weyl $m$-function}\label{m-function}

Let $v_\beta$ be the solution of (\ref{diff}) on $[-S,S]$ obeying the transfer condition (\ref{tx-condition}) and satisfying the terminal conditions $v_\beta(S,\lambda )= \sin \beta $ and $v_\beta'(S,\lambda ) = \cos \beta $.
To define the Titchmarsh-Weyl $m$-function we use the approach given in \cite{BBW}, i.e.
the $m$-function of (\ref{diff}) on $[-S,S]$ for boundary conditions (\ref{alpha}), (\ref{beta})
and the transfer condition (\ref{tx-condition}) 
is that value of $m_{\alpha,\beta}$ for which 
\begin{equation}\label{mfn}
\psi_{\alpha,\beta} :=u_\alpha + m_{\alpha,\beta}w_\alpha
\end{equation}
obeys the  terminal condition (\ref{beta}). Here $u_\alpha, w_\alpha$ are solutions of (\ref{diff}) obeying the initial condition
\begin{equation}
 W_\alpha(-S,\lambda ) = H_{\alpha}= 
\left(\begin{array}{cc}
 \cos \alpha & \sin \alpha \\
-\sin \alpha & \cos \alpha
  \end{array} \right)
\label{halpha}	
\end{equation}	
where
\begin{equation}\label{X}
W_\alpha(x, \lambda ) := \left[\begin{array}{cc}
u_\alpha(x,\lambda ) & w_\alpha(x,\lambda )\\
u_\alpha'(x,\lambda ) & w_\alpha'(x,\lambda )
\end{array}\right].
\end{equation}
The entries of $W_\alpha(x,\lambda )$ are entire functions of $\lambda$ and the determinant  of $W_\alpha(x,\lambda)$ is the Wronskian of $u_\alpha$ and $w_\alpha$, which is $1$ for all $x$ and $\lambda.$

Let 
\begin{eqnarray}
\Delta_{\alpha,\beta} (\lambda) := \makebox{Wron}[w_\alpha,v_\beta] = w_\alpha v_\beta' - v_\beta w_\alpha' = v_\beta'(-S,\lambda)\sin \alpha -v_\beta(-S,\lambda)\cos \alpha .\label{delta-baw}
\end{eqnarray}
The function $\Delta_{\alpha,\beta} (\lambda )$ is entire in $\lambda $ and the zeros of $\Delta_{\alpha,\beta} (\lambda )$ are the eigenvalues of (\ref{diff}) with  boundary conditions (\ref{alpha}) and (\ref{beta}) and the
transfer condition (\ref{tx-condition}).
For $\lambda$ not an eigenvalue of (\ref{diff}) -- (\ref{beta}), $v_\beta$ and $\psi_{\alpha,\beta}$ are linearly dependent, say $v_\beta=k \psi_{\alpha,\beta}$, giving 
\begin{eqnarray*}
v_\beta(-S,\lambda)&=& -k\cos\alpha+ km_{\alpha,\beta}(\lambda)\sin\alpha,\\
v_\beta'(-S,\lambda)&=&k\sin\alpha + km_{\alpha,\beta}(\lambda)\cos\alpha.
\end{eqnarray*}
It thus follows that $k=v_\beta'(-S,\lambda)\sin\alpha-v_\beta(-S,\lambda)\cos\alpha=\Delta_{\alpha,\beta}(\lambda)$ and
\[\psi_{\alpha,\beta} (x,\lambda ) 
= \frac{v_\beta(x,\lambda )}{\Delta_{\alpha,\beta}\lambda )}.\]
Forming the linear combinations
$v_\beta(-S,\lambda)\sin\alpha+v_\beta'(-S,\lambda)\cos\alpha=km_{\alpha,\beta}(\lambda)$ we have
\begin{eqnarray}
m_{\alpha,\beta}(\lambda)
&=&\frac{v_{\beta}(-S,\lambda)\sin \alpha +v_{\beta}'(-S,\lambda)\cos \alpha}{\Delta_{\alpha,\beta}(\lambda)}. \label{mfngen}
\end{eqnarray}


\begin{lem}\label{mfnasymp}
The Titchmarsh-Weyl m-function for $\alpha=0$ and $\beta =\pi$, $m_{0,\pi}$, 
 has asymptotic approximation
\begin{equation*}
 m_{0,\pi}(\lambda ) = -i\sqrt{\lambda} + O(1)\quad\mbox{as}\quad \lambda \rightarrow -\infty.
\end{equation*}
\end{lem}

\proof
From (\ref{delta-baw}) and (\ref{mfngen}), for $\alpha=0$,
\begin{equation}\label{mfnv}
m_{0,\beta}(\lambda) = - \frac{v_\beta'(-S,\lambda)}{v_\beta(-S, \lambda)}.
\end{equation}
Denote the entries of the  transfer matrix $M$  by $m_{ij}, i,j = 1,2.$
Note that, as $\det M =1$ and $M$ has real entries not both $m_{12}$ and $m_{11} + m_{22}$ can simultaneously be zero.
Let $\sqrt{\lambda}=ik, k>0$, then
from \cite[Appendix]{cenw2}, as $k\to \infty$,
\begin{eqnarray}
v_\pi(-S,-k^2)&=&-\frac{km_{12}+m_{11}+m_{22}}{4k^2}e^{2Sk}[k+O(1)],\label{v-beta-zero-1}\\
v_\pi'(-S,-k^2)&=&
\frac{km_{12}+m_{11}+m_{22}}{4k}e^{2Sk}[k+O(1)].\label{v-beta-zero-2}
\end{eqnarray}
Therefore, as $k \rightarrow \infty$,
\begin{eqnarray*}
-\frac{v_\pi'(-S,-k^2)}{v_\pi(-S,-k^2)}=k+O(1)=-i\sqrt{\lambda}+O(1).
 \qed
\end{eqnarray*}


\begin{lem}\label{mfnasymp1}
The Titchmarsh-Weyl m-function for $\alpha=0$ and $\beta\in (0,\pi)$, $m_{0,\beta}$, 
 has asymptotic approximation
\begin{equation*}
 m_{0,\beta}(\lambda ) = -i\sqrt{\lambda} + O(1)\quad\mbox{as}\quad \lambda \rightarrow -\infty.
\end{equation*}
\end{lem}

\proof
We proceed from (\ref{mfnv}), again let  $\sqrt{\lambda}=ik, k>0$, then
from \cite[Appendix]{cenw2}, as $k\to \infty$, for $\beta\in (0,\pi)$,
\begin{eqnarray}
v_\beta(0^+,-k^2) &=& -\frac{\sin \beta}{2}e^{kS}\left( 1 +O\left(\frac{1}{k} \right)\right), \label{214}\\
v_\beta'(0^+,-k^2) &=& k\frac{\sin \beta}{2}e^{kS}\left(1 +O\left(\frac{1}{k} \right) \right). \label{215}
\end{eqnarray}
Since $v_\beta(x,\lambda)$ satisfies transfer condition (\ref{tx-condition}) we have 
\begin{eqnarray}
v_\beta(0^-,-k^2) &=& -(km_{12}+m_{22})\frac{\sin \beta}{2}e^{kS}\left(1 +O\left(\frac{1}{k} \right) \right), \\
v_\beta'(0^-,-k^2) &=& (km_{11}+m_{21})\frac{\sin \beta}{2}e^{kS}\left(1 +O\left(\frac{1}{k} \right) \right).
\end{eqnarray}
On $(-S,0)$,  
$$v_\beta(x,-k^2)=v_\beta(0^-,-k^2)h_0(x)+v_\beta'(0^-,-k^2)h_1(x), $$
where $h_0, h_1$ are solutions of (\ref{diff}) with 
$h_0(0^-)=1=h_1'(0^-)$ and $h_0'(0^-)=0=h_1(0^-)$.
Here, as $\lambda=-k^2\to -\infty$,
$h_0(x)=\frac{e^{-kx}}{2}\left(1+O\left(\frac{1}{k} \right)\right)$
and $h_1(x)=-\frac{e^{-kx}}{2k}\left(1+O\left(\frac{1}{k} \right)\right)$ with
$h_0'(x)=-k\frac{e^{-kx}}{2}\left(1+O\left(\frac{1}{k} \right)\right)$
and $h_1'(x)=\frac{e^{-kx}}{2}\left(1+O\left(\frac{1}{k} \right)\right)$.
Thus
\begin{eqnarray}
v_\beta(-S,-k^2)&=&
-(k^2m_{12}+k(m_{11}+m_{22})+m_{21})\frac{\sin \beta}{4k}e^{2kS}
\left(1+O\left(\frac{1}{k} \right)\right), \label{v-beta-non-zero-1}\\
v_\beta'(-S,-k^2)&=&
k(k^2m_{12}+k(m_{11}+m_{22})+m_{21})\frac{\sin \beta}{4k}e^{2kS}\left(1+O\left(\frac{1}{k} \right)\right).\label{v-beta-non-zero-2}
\end{eqnarray}
Thus $m_{0,\beta}(-k^2)=k + O(1).$
\qed


\begin{lem}\label{mfnasymp2}
The Titchmarsh-Weyl m-function given in (\ref{mfn}) has the following asymptotic approximation for $\alpha,\beta\in(0,\pi)$  as $\lambda \rightarrow -\infty$
\begin{equation}\label{mfnasympeqn2}
 m_{\alpha,\beta}(\lambda ) = \cot \alpha + O\left(\frac{1}{\sqrt{\lambda}}\right).
\end{equation}
\end{lem}

\proof
Let $\lambda=-k^2$ and $k>0$ for $k\to\infty$ we have $v_\beta(-S,-k^2)$ and $v'_\beta(-S,-k^2)$
as given by (\ref{v-beta-non-zero-1})-(\ref{v-beta-non-zero-2}).
Direct computation give
$$ \Delta_{\alpha,\beta}(\lambda)=
\frac{k^2m_{12}+k(m_{11}+m_{22})}{4} \sin \alpha \sin \beta  e^{2kS}\left(1+O\left(\frac{1}{k}\right)\right),$$
and
\begin{eqnarray*}
v_{\beta}(-S,\lambda)\sin \alpha +v_{\beta}'(-S,\lambda)\cos \alpha
=\frac{k^2m_{12}+k(m_{11}+m_{22})}{4} \cos \alpha \sin \beta  e^{2kS}\left(1+O\left(\frac{1}{k}\right)\right).
\end{eqnarray*}
The result thus follows from (\ref{mfngen}).
\qed


\begin{lem}\label{mfnasymp3}
The Titchmarsh-Weyl m-function given in (\ref{mfn}) has the following asymptotic approximation for $\alpha\in (0,\pi)$ and $\beta = \pi$ as $\lambda \rightarrow -\infty$
\begin{equation}\label{mfnasympeqn3}
 m_{\alpha,\pi}(\lambda ) = \cot\alpha + O\left(\frac{1}{\sqrt{\lambda}}\right).
\end{equation}
\end{lem}

\proof
Let $\lambda=-k^2$ for $k>0$, then for $k\to\infty$ from (\ref{v-beta-zero-2}) and (\ref{v-beta-zero-2}) we have
$$\Delta_{\alpha,\pi} (\lambda)
= \sin\alpha\left(km_{12}+m_{11}+m_{22})(k+\cot \alpha) \right)\frac{e^{2Sk}}{4k^2}[k+O(1)]$$
and
$$v_{\beta}(-S,\lambda)\sin \alpha +v_{\beta}'(-S,\lambda)\cos \alpha
=\sin\alpha(km_{12}+m_{11}+m_{22})(k\cot \alpha-1)\frac{e^{2Sk}}{4k^2}[k+O(1)].$$
Hence
\begin{equation*}
 m_{\alpha,\pi}(\lambda ) = 
\frac{k\cot \alpha-1}{k+\cot \alpha}\left(1+O\left(\frac{1}{k}\right)\right). \qed
\end{equation*}


Combining the four lemmas above we can write the following.
\begin{thm} \label{cor}
The Titchmarsh-Weyl m-function given in (\ref{mfn}) has the following asymptotic behaviour as  $\lambda \rightarrow -\infty$,
\begin{equation}\label{mfnasympeqn4}
 m_{\alpha,\beta}(\lambda ) = \left\{ 
\begin{array}{ll}
-  i\sqrt{\lambda} + O(1), & \alpha =0, \\
\cot \alpha +O\left(\frac{1}{\sqrt{\lambda}} \right),&  \alpha \neq 0.
\end{array}
\right.
\end{equation}
\end{thm}
Note that this result is consistent with that obtained in \cite{BBW}.

\newsection{Inverse problems on a compact interval}

In this section we consider various inverse spectral problems for the differential equation
(\ref{diff}) on the interval $[-S,S]$ with transmission condition (\ref{tx-condition}) at $0$
and boundary conditions (\ref{alpha}) and (\ref{beta}).

\begin{thm}\label{m2spec}
Let $\lambda_0<\lambda_1<\lambda_2<\dots$ be the eigenvalues of (\ref{diff}) and (\ref{tx-condition}) with boundary conditions (\ref{alpha}) and (\ref{beta}).
Let $\mu_0<\mu_1<\mu_2<\dots$ be the eigenvalues of (\ref{diff}) and (\ref{tx-condition}) with boundary conditions (\ref{alpha}) and (\ref{beta}), where $\alpha$ in (\ref{alpha}) has been replaced by $\alpha'\ne \alpha, \alpha'\in [0,\pi)$.
The Titchmarsh-Weyl m-function, $m_{\alpha,\beta}$, can be uniquely reconstructed from $(\lambda_n), (\mu_n), \alpha$ and $\alpha'$.
\end{thm}

\proof 
The functions $\Delta_{\alpha,\beta}(\lambda)$ and $\Delta_{\alpha',\beta}(\lambda)$ are entire of order $\frac12$ in $\lambda$. The zeros of the former are $(\lambda_n)$ while the zeros of the latter are $(\mu _n)$. Thus using the Hadamard product theorem we can write
$$ \Delta_{\alpha,\beta}(\lambda)=C_{\alpha,\beta}\prod_{n=0}^{\infty}\left(1-\frac{\lambda}{\lambda_n} \right) $$
and 
$$ \Delta_{\alpha',\beta}(\lambda)=C_{\alpha',\beta}\prod_{n=0}^{\infty}\left(1-\frac{\lambda}{\mu_n} \right). $$
Note that if for some $n$, $\lambda_n=0$ then we consider $\frac{\Delta_{\alpha,\beta}(\lambda)}{\lambda}$ instead. Similarly if $\mu_n =0$ for some $n$ consider $\frac{\Delta_{\alpha',\beta}(\lambda)}{\lambda}$. For convenience we denote

Now 
$$\cot(\alpha-\alpha')\left[\begin{array}{c} \sin\alpha\\ \cos\alpha\end{array}\right]
-\mbox{cosec}(\alpha-\alpha') \left[\begin{array}{c} \sin\alpha'\\ \cos\alpha'\end{array}\right]
=\left[\begin{array}{c} \cos\alpha\\ -\sin\alpha\end{array}\right]=\left[\begin{array}{c} \sin(\frac{\pi}{2}+\alpha)\\ \cos(\frac{\pi}{2}+\alpha)\end{array}\right],$$ 
so
$$\cot(\alpha-\alpha')\Delta_{\alpha,\beta}
-\mbox{cosec}(\alpha-\alpha')\Delta_{\alpha',\beta}
=\Delta_{\alpha+\frac{\pi}{2},\beta}.$$
Hence
$$m_{\alpha,\beta}
=\frac{\Delta_{\alpha+\frac{\pi}{2},\beta}}{\Delta_{\alpha,\beta}}
=\cot(\alpha-\alpha')
-\mbox{cosec}(\alpha-\alpha')\frac{\Delta_{\alpha',\beta}}{\Delta_{\alpha,\beta}}$$
giving
$$m_{\alpha,\beta}(\lambda) =\cot(\alpha-\alpha')
-C\mbox{cosec}(\alpha-\alpha')
\frac{ \prod_{n=0}^{\infty} \left(1-\frac{\lambda }{\mu _n}\right)}{\prod_{n=0}^{\infty} \left(1-\frac{\lambda }{\lambda _n}\right)}$$
where $C=C_{\alpha',\beta}/C_{\alpha,\beta}$.  In particular if $\alpha,\alpha'\ne 0$ then
$$\lim_{\lambda\to -\infty}\frac{\sin\alpha'}{\sin\alpha}
\frac{\prod_{n=0}^{\infty} \left(1-\frac{\lambda }{\lambda _n}\right)}
{ \prod_{n=0}^{\infty} \left(1-\frac{\lambda }{\mu _n}\right)}
=C$$
while if $\alpha=0$ then $\alpha'\in (0,\pi)$ and
$$\lim_{\lambda\to -\infty}-i\sqrt{\lambda}\sin\alpha'
\frac{\prod_{n=0}^{\infty} \left(1-\frac{\lambda }{\lambda _n}\right)}
{ \prod_{n=0}^{\infty} \left(1-\frac{\lambda }{\mu _n}\right)}
=C,$$
from Corollary \ref{cor}. 

Finally, for $\alpha'=0$ and $\alpha\in (0,\pi)$ to solve for $C$ we need to consider the asymptotic
behaviour of $\Delta_{0,\beta}(\lambda)/\Delta_{\alpha,\beta}(\lambda)$ as $\lambda\to -\infty$.
From the previously computed asymptotic approximations, as $k\to \infty$,
$$\Delta_{0,\beta}(-k^2)=\left\{\begin{array}{ll} 
\frac{km_{12}+m_{11}+m_{22}}{4k^2}e^{2Sk}[k+O(1)],&\beta=\pi\\
(km_{12}+m_{11}+m_{22})\frac{\sin \beta}{4k}e^{2kS}
[k+O(1)],&\beta\in (0,\pi)
\end{array}\right.$$
and
$$ \Delta_{\alpha,\beta}(-k^2)=\left\{\begin{array}{ll}
\sin\alpha\left(km_{12}+m_{11}+m_{22})(k+\cot \alpha) \right)\frac{e^{2Sk}}{4k^2}[k+O(1)],&
\beta=\pi\\
\frac{km_{12}+m_{11}+m_{22}}{4} \sin \alpha \sin \beta  e^{2kS}[k+O(1)],&\beta\in (0,\pi)\\
\end{array}\right.$$
giving
$$\frac{\Delta_{0,\beta}(\lambda)}{\Delta_{\alpha,\beta}(\lambda)}
=\frac{1}{-i\sqrt{\lambda}\sin\alpha}\left(1+O\left(\frac{1}{\sqrt{\lambda}}\right)\right),$$
as $\lambda\to -\infty$ and
$$\lim_{\lambda\to -\infty} \frac{1}{-i\sqrt{\lambda}\sin\alpha} 
\frac
{\prod_{n=0}^{\infty} \left(1-\frac{\lambda }{\lambda _n}\right)}
{ \prod_{n=0}^{\infty} \left(1-\frac{\lambda }{\mu _n}\right)}
=C.
\qed$$


We now show that the Titchmarsh-Weyl $m$-function can be determined from one spectrum together with the associated norming constants.  Here we recall that if $\lambda_n$ is an eigenvalue of 
(\ref{diff})-(\ref{beta}) then the norming constant associated with the eigenvalue $\lambda_n$ is
$$a_n=\int_{-S}^S w_\alpha^2(x,\lambda_n)\,dx.$$

\begin{thm}\label{mfn1specnorm}
For  (\ref{diff}) on $[-S,S]$ with transmission condition (\ref{tx-condition}) at $x=0$ and
boundary conditions (\ref{alpha})-(\ref{beta}) 
the eigenvalues $\lambda_0<\lambda_1<\lambda_2<\dots$  together with the corresponding norming constants $(a_n)$, uniquely determine the Titchmarsh-Weyl m-function, $m_{\alpha ,\beta }$ for $\alpha \ne 0$.
\end{thm}

\proof 
Since the eigenvalues are simple, for each $n=0,1,2,\dots,$ there is a non-zero $k_n$ so that $k_nw_\alpha(x,\lambda _n) = v_\beta(x,\lambda _n)$. 
and
$k_nw_\alpha'(x,\lambda _n) = v_\beta'(x,\lambda _n)$.
Hence
$$k_n=k_n \mbox{Wron}[u_\alpha,w_\alpha](x,\lambda _n)=\mbox{Wron}[u_\alpha,v_\beta](x,\lambda _n)=\Delta_{\alpha+\frac{\pi}{2},\beta}(\lambda_n).$$

Since $w_\alpha(x,\lambda_n )$ and $v_\beta(x,\lambda)$ both obey (\ref{beta}) we have
$\makebox{Wron}[v_\beta(x,\lambda),w_\alpha(x,\lambda _n)](S)=0$ and thus
\begin{eqnarray*}
 \int_{-S}^S (\lambda -\lambda _n)v_\beta(x,\lambda )w_\alpha(x,\lambda _n)\,dx &=& \int_{-S}^S [w_\alpha''(x,\lambda_n )v_\beta(x,\lambda) -w_\alpha(x,\lambda _n)v_\beta''(x,\lambda )]\, dx\\
 &=& \left[\makebox{Wron}[v_\beta(x,\lambda),w_\alpha(x,\lambda _n)]\right]_{-S}^S\\
 &=& \Delta_{\alpha\beta}  (\lambda ).
\end{eqnarray*}
So, as $\Delta_{\alpha\beta} (\lambda_n)=0$,
\[\frac{\Delta_{\alpha\beta} (\lambda )-\Delta_{\alpha\beta} (\lambda_n)}{\lambda -\lambda _n} = \int_{-S}^S v(x,\lambda )w_2(x,\lambda _n)\,dx.\]
Thus taking $\lambda \rightarrow \lambda _n$, we have
\[\Delta_{\alpha\beta} '(\lambda _n) = \int_{-S}^S v_\beta(x,\lambda _n)w_\alpha(x,\lambda _n)\,dx =k_na_n=\Delta_{\alpha+\frac{\pi}{2},\beta}(\lambda_n)a_n.\]
Combining the above with the relation
$$m_{\alpha,\beta}
=\frac{\Delta_{\alpha+\frac{\pi}{2},\beta}}{\Delta_{\alpha,\beta}}
$$
(where $\Delta_{\alpha+\frac{\pi}{2},\beta}:=\Delta_{(\alpha+\frac{\pi}{2}){\rm mod}\pi,\beta}$)
gives that
\[\makebox{Res}_{\lambda =\lambda _n} m(\lambda )= \frac{1}{a_n}.\]

Without loss of generality we may assume that $0$ is not an eigenvalue (if it is shift the potential and the eigenparameter so that in the shifted eiegnparameter $0$ is not an eigenvalue).  For $\alpha\in (0,\pi)$ and
$\beta \in (0,\pi]$, using the asymptotic estimates from the Appendix, it follows that for large $n\in\N$,
$|m_{\alpha,\beta}(\lambda)|\le  |\cot\alpha|+1$ for on circle $C_n$ with $|\lambda|=\left(n+\frac{1}{4}\right)\pi$.  Hence the Mittag-Leffler expansion theorem can be applied to give
\begin{eqnarray}
m_{\alpha,\beta}(\lambda)=m_{\alpha,\beta}(0)+\sum_{n=0}^\infty \frac{1}{a_n}\left(\frac{1}{\lambda-\lambda_n}+\frac{1}{\lambda_n}  \right).\label{m-2015-july}
\end{eqnarray}  
Applying Theorem \ref{cor} to (\ref{m-2015-july}) enables one to solve for the constant
$m_{\alpha,\beta}(0)$. \qed

{\bf Remark} For $\alpha=0$, let $\mu<\min\{0,\lambda_0\}$. 
Using the asymptotic estimates given in the Appendix it can be seen that the function 
$g(\lambda):=\frac{m_{0,\beta}(\lambda)}{\lambda-\mu}$ satisfies the conditions of the Mittag-Leffler expansion theorem.
The poles of $g(\lambda)$ are $\mu<\lambda_0<\lambda_1<\dots$ and a direct computation gives that the
residues of $g$ at these poles are $m_{0,\beta}(\mu), \frac{1}{a_0(\lambda_0-\mu)}, \frac{1}{a_1(\lambda_1-\mu)}, \dots$.  Hence
\begin{eqnarray}\label{g-2015}
 g(\lambda)= -\frac{m_{0,\beta}(0)}{\mu}
+m_{0,\beta}(\mu)\left(\frac{1}{\lambda-\mu}+\frac{1}{\mu}\right)
+\sum_{n=0}^\infty \frac{1}{a_n(\lambda_n-\mu)}  
   \left(\frac{1}{\lambda-\lambda_n}+\frac{1}{\lambda_n}  \right).
\end{eqnarray}
Multiplying by $\lambda-\mu$ we obtain
\begin{eqnarray}\label{m-0-2015}
 m_{0,\beta}(\lambda)=
m_{0,\beta}(0)
+(m_{0,\beta}(\mu)-  m_{0,\beta}(0) )\frac{\lambda}{\mu}
+\sum_{n=0}^\infty
   \frac{\lambda}{a_n\lambda_n}\left(
 \frac{1}{\lambda_n-\mu}+   \frac{1}{\lambda-\lambda_n}\right).
\end{eqnarray}
As $m_{0,\beta}(\mu)$ is analytic at $\mu=0$ we can take
$\mu\to 0$ in the above, giving
\begin{eqnarray}\label{m-1-2015}
 m_{0,\beta}(\lambda)=
m_{0,\beta}(0)
+\lambda\left[m_{0,\beta}'(0)
+\sum_{n=0}^\infty
   \frac{1}{a_n\lambda_n}\left(
 \frac{1}{\lambda_n}+   \frac{1}{\lambda-\lambda_n}\right)\right].
\end{eqnarray}

Dividing by $\lambda$ and taking $\lambda\to -\infty$ in (\ref{m-1-2015}), from Theorem \ref{cor}
we have 
\begin{equation}\label{anln} 
m_{0,\beta}'(0)=-
\sum_{n=0}^\infty
   \frac{1}{a_n\lambda_n^2}.
   \end{equation}
Using the asymptotics in the Appendix for $v(x,\lambda _n)$ together with the fact that $\lambda _n = O(n^2)$ one can show that $\displaystyle{\frac{1}{a_n\lambda _n^2} = O\left(n^{-2}\right)}$. Therefore the series in (\ref{anln}) converges and thus
\begin{eqnarray}\label{m-2-2015}
 m_{0,\beta}(\lambda)=
m_{0,\beta}(0)
+\lambda
\sum_{n=0}^\infty
   \frac{1}{a_n\lambda_n(\lambda-\lambda_n)}\,\,.
\end{eqnarray}
Note that if the constant term $C$, for $\alpha =0$, in (\ref{mfnasympeqn4}) were known i.e.
$$m_{0,\beta }(\lambda ) = -i\sqrt{\lambda } + C + O\left(\frac{1}{\sqrt{\lambda }}\right),$$
then $m_{0,\beta }(0) $ in (\ref{m-2-2015}) is uniquely given by 
$$m_{0,\beta }(0) = \lim_{\lambda \rightarrow -\infty} \left(-i\sqrt{\lambda } - \lambda
\sum_{n=0}^\infty
   \frac{1}{a_n\lambda_n(\lambda-\lambda_n)}\right) + C.$$
Whether or not $C$ and $\,\, \displaystyle{\sum_{n=0}^\infty
   \frac{1}{a_n\lambda_n(\lambda-\lambda_n)}}\,\,$ are practically accessible remains an open question.

\newsection{Scattering problem on the line}

In this section we consider the inverse scattering problem on the real line i.e. for the differential equation
\begin{equation}\label{diffapp}
 \ell y= {\zeta}^2 y,\quad \mbox{where}\quad \ell y:=-\frac{d^2y}{dx^2} + q(x)y,\quad \mbox{on } (-\infty,0)\cup (0,\infty),
\end{equation}
with point transfer condition
\begin{equation}\label{tx-conditionapp}
  \left[\begin{array}{c}y(0^+)\\ y'(0^+)\end{array}\right]=M
   \left[\begin{array}{c}y(0^-)\\ y'(0^-)\end{array}\right].
\end{equation}
Here the entries of $M$ are taken to be real, $q \in L^2(\R)$ is assumed to be real valued and have compact 
essential support, say ${\rm ess\, supp}(q) \subset [-S,S]$ for some $S>0$.
The scattering problem can now be treated as two half-line problems 
interacting via a matrix transfer condition (\ref{tx-condition}) at the origin. 

\begin{defs}\label{def_Jost}
\cite[p.297]{chadan} The Jost solutions $f_{+,M}(x,\zeta)$ and $f_{-,M}(x,\zeta)$ are the solutions
 of (\ref{diffapp}) and (\ref{tx-conditionapp}) with
\begin{equation}
\lim_{x \to \infty} e^{-i\zeta x} f_{+,M}(x,\zeta) = 1
 = \lim_{x \to -\infty} e^{i\zeta x} f_{-,M}(x,\zeta). 
\end{equation}
\end{defs}

The Jost solutions  $f_{+,M}(x,\zeta)$ and $f_{-,M}(x,\zeta)$ to 
(\ref{diffapp}), (\ref{tx-conditionapp}), can be expressed in terms of  the
 classical Jost solutions $f_+(x,\zeta)$ and $f_-(x,\zeta)$ (i.e. when $M=I$)
 by
\begin{equation}\label{F+}
f_{+,M}(x,\zeta) := \left\{\begin{array}{c}f_+(x,\zeta), \quad x > 0 \\
h_1(x,\zeta), \quad x < 0 \end{array}\right.,
\end{equation}
\begin{equation}\label{F-}
f_{-,M}(x,\zeta) := \left\{\begin{array}{c}f_-(x,\zeta), \quad x < 0 \\
h_2(x,\zeta), \quad x > 0 \end{array}\right.,
\end{equation}
where $h_1(x,\zeta)$ and $h_2(x,\zeta)$ are solutions of (\ref{diffapp}) on $(-\infty, 0)$ and $(0, \infty)$ respectively obeying 
\begin{eqnarray*}
\left(\begin{array}{c}h_1(0^-, \zeta) \\ h'_1(0^-,\zeta)\end{array}\right)&=& M^{-1}\left(\begin{array}{c}f_+(0^+, \zeta) \\ f'_+(0^+,\zeta)\end{array}\right),\\
\left(\begin{array}{c}h_2(0^+, \zeta) \\ h'_2(0^+,\zeta)\end{array}\right)&=&M\left(\begin{array}{c}f_-(0^-, \zeta) \\ f'_-(0^-,\zeta)\end{array}\right).
\end{eqnarray*}

For $M=I$ the existence and asymptotic behaviour of the Jost solutions have been well studied, 
see for example \cite{Freiling, March}. In particular
\begin{eqnarray}\label{4_1f}
 f_+(x,\zeta) &=& e^{i{\zeta}x} + O\left(\frac{C(x)\rho(x)e^{-{\eta}x}}{1 +
|\zeta|}\right), \\
\label{4_3}
f_-(x,\zeta) &=& e^{-i{\zeta}x} + O\left(\frac{C(-x)\tilde{\rho}(x)e^{{\eta}x}}{1 +
|\zeta|}\right), 
\end{eqnarray}
as $|x|+|\zeta| \rightarrow \infty$, where $\eta = \Im (\zeta)$. Here $C(x)$ is a non-negative, non-increasing function of $x$ and
\begin{equation}
 \rho(x) = \int_x^{\infty} (1+|\tau |)|q(\tau)|\,d \tau,
\qquad \tilde{\rho}(x) = \int_{-\infty}^x (1+|\tau |)|q(\tau)|\,d \tau .
\end{equation}
For $\xi \in \R$, see \cite[Sections 2 and 4]{cenw1}, 
the conjugate Jost solutions take the form
\begin{equation}\label{big_F_conj}
 \overline{f}_{+,M}(x,\xi):= 
   \left\{\begin{array}{cc}\overline{f}_+(x,\xi) = f_+(x,-\xi),\quad & x>0 \\
    \overline{h}_1(x,\xi) = h_1(x,-\xi),\quad & x<0 \end{array}\right.
\end{equation}
which obeys the transfer condition at $x=0$.
The solutions $f_{+,M}(x,\xi)$ and 
$\overline{f}_{+,M}(x,\xi)$ are independent for  $\xi\in\R\backslash\{0\}$ and
span the solution space of (\ref{diffapp}), with (\ref{tx-conditionapp}).  Hence
there are unique coefficients $A(\xi)$ and $B(\xi)$ so that
\begin{equation}\label{Big_A_and_B}
f_{-,M}(x,\xi) = A(\xi)\overline{f}_{+,M}(x,\xi) + B(\xi)f_{+,M}(x,\xi).
\end{equation}
Here $A(\xi)$ and $B(\xi)$ are independent of whether $x>0$ or $x<0$, and
they satisfy the equality
$|A(\xi)|^2 - |B(\xi)|^2 = 1$ for $\xi\in\R\backslash\{0\}$.
The reflection coefficient is defined as
\begin{eqnarray*}
  R(\xi) = \frac{B(\xi)}{A(\xi)}, \mbox{ for } \xi \in \R\backslash\{0\}.
\end{eqnarray*}

The next two results are a generalisation of \cite[Lemma 4.1, Theorem 4.2]{cenw2}. The proofs follow in exactly the same manner, we will point out the main differences and provide the necessary asymptotics.

\begin{lem} \label{4.1new}
Let $\mbox{ess supp}(q)\in [-S,S]$ for some $S>0.$ Given the scattering data $\{R(\xi), \eta_1, \dots, \eta_N \}$, the matrix $W(S,\zeta)$ is uniquely determined. Here, $W(x,\zeta)$ is as given in (\ref{X}) with $W(-S,\zeta)$ defined in (\ref{halpha}).
\end{lem}

\proof
Similarly to the proof of \cite[Lemma 4.1]{cenw2} we now obtain 
\begin{equation}
[w_1(S,\xi)\; w_2(S,\xi)] = [e^{-i\xi S}\; e^{i\xi S}] 
\begin{pmatrix}
A(\xi) & \overline{B}(\xi) \\
B(\xi) & \overline{A}(\xi)	
\end{pmatrix}
\begin{pmatrix}
	\frac{e^{-i\xi S}}{2} & -\frac{e^{-i\xi S}}{2i\xi} \\
	\frac{e^{i\xi S}}{2} & \frac{e^{i\xi S}}{2i \xi}
\end{pmatrix}
H_{\alpha}.
\label{w1w2}
\end{equation}
Therefore as in \cite{cenw2} we can find $w_1(S,\xi)$ and $w_2(S,\xi)$ and since $w_1$ and $w_2$ are entire, by analyticity we can extend them to $w_1(S,\zeta)$ and $w_2(S,\zeta).$ \qed

\begin{thm}\label{potqnew}
Given the Titchmarsh-Weyl m-function, $m$, to (\ref{diff}) on $[-S,S]$ with
 boundary conditions (\ref{alpha}) and (\ref{beta}) and the transfer condition (\ref{tx-condition}) and $\tilde{m}$,
 the Titchmarsh-Weyl m-function for the same problem but with the potential $q$ replaced by
 $\tilde{q}$. If $m = \tilde{m}$ then $q= \tilde{q}$.
\end{thm}

\proof 
The proof follows identically to that given in \cite[Theorem 4.2]{cenw2}. It should be noted that only the asymptotics for $v$ and $w_2$ used in Theorems 5.3 and 5.4 change. The new required asymptotics for $v$ are given by (\ref{214}),(\ref{215}) and (\ref{v1})- (\ref{v4}). The asymptotics for $w_2$ are as follows:

For $-S \leq x <0 $ we get
\begin{eqnarray}
w_2 (x,\lambda) &=&\sin \alpha  \cos \sqrt{\lambda}(x+S) + 
O \left( \frac{e^{|\Im \sqrt{\lambda}|(x+S)}}{\sqrt{\lambda}}\right),
\label{w21} \\
w_2'(x,\lambda) &=& -\sqrt{\lambda} \sin \alpha \sin \sqrt{\lambda}(x+S)+O(e^{|\Im \sqrt{\lambda}|(x+S)}). \label{w22}
\end{eqnarray}
For $S \geq x >0$ we have if $m_{12}\neq 0$ then 
\begin{eqnarray}
w_2(x,\lambda) &=& -m_{12} \sqrt{\lambda} \sin \alpha \sin \sqrt{\lambda} S \cos \sqrt{\lambda}x + O (e^{|\Im \sqrt{\lambda}|(x+S)}), \label{w23} \\ 
w_2(x,\lambda) &=& m_{12} {\lambda} \sin \alpha \sin \sqrt{\lambda} S \sin \sqrt{\lambda}x + O (\sqrt{\lambda} e^{|\Im \sqrt{\lambda}|(x+S)}), \label{w24}
\end{eqnarray}
and if $m_{12}=0$ then
\begin{eqnarray}
w_2(x,\lambda) &=& m_{11} \sin \alpha \cos \sqrt{\lambda} x \cos \sqrt{\lambda}S - m_{22} \sin \alpha \sin \sqrt{\lambda} x \sin \sqrt{\lambda} S \nonumber \\ 
& & + O  \left( \frac{e^{|\Im \sqrt{\lambda}|(x+S)}}{\sqrt{\lambda}}\right),
\label{w25} \\
w_2'(x,\lambda) &=& \sqrt{\lambda} \sin \alpha(-m_{11} \sin \sqrt{\lambda} x \cos \sqrt{\lambda}S - m_{22} \cos \sqrt{\lambda} x \sin \sqrt{\lambda} S ) \nonumber \\
& & + O \left( {e^{|\Im \sqrt{\lambda}|(x+S)}}\right).
\label{w26}
\end{eqnarray}
\hfill \qed

Combining Theorem \ref{m2spec} with Theorem \ref{potqnew} gives the following.

\begin{cor}\label{cor1}
Given two spectra where one spectrum comes from the problem (\ref{diff}), (\ref{tx-condition}) with boundary conditions (\ref{alpha}) and (\ref{beta}) and the other spectrum comes from the same problem but with a different $\beta \in (0,\pi]$ at the terminal end point, the potential $q$ can be uniquely determined on $[-S,S]$.
\end{cor}

Also combining Theorem \ref{mfn1specnorm} with Theorem \ref{potqnew} gives the result below.

\begin{cor}\label{cor1.1}
Given the spectrum of the Neumann-Neumann problem and corresponding norming constants for equation (\ref{diff}) on $[-S,S]$ and transfer matrix, $M$, at $x=0$, the potential $q$ may be uniquely determined on $[-S,S]$.
\end{cor}

\begin{thm}\label{fin}
Given the scattering data $\{R(\xi), \eta_1, \dots, \eta_N \}$, i.e. the reflection coefficient and the eigenvalues of the problem (\ref{diffapp}), (\ref{tx-conditionapp}), the spectra of the Neumann-Neumann and the Neumann-Dirichlet problems on $[-S,S]$, with potential $q|_{[-S,S]}$ and transfer matrix, $M$, as in the original scattering problem on the line, can be obtained. Moreover, the norming constants for the Neumann-Neumann problem on $[-S,S]$ can be determined.
\end{thm}

\proof
Since $w_1(x,\zeta)$ and $w_2(x,\zeta)$ build a fundamental system and we know their values at $x=-S$ and from Lemma 
\ref{4.1new} we can determined $w_1(S,\zeta)$ and $w_2(S,\zeta),$ we can obain (from the scattering data) any spectrum, in particular the Neumann-Neumann spectrum, say $\lambda_n.$ Thus we can construct $\Delta(\lambda)$ for
$\alpha = \frac{\pi}{2}$ and therefore we know $\Delta'(\lambda)$ for Neumann-Neumann boundary conditions. 
In addition from Lemma \ref{4.1new} we know $w_2'(S,\lambda)$ and 
$$ 
\mbox{Wron}[w_1,v](S,\lambda) = -w_1'(S,\lambda) = \mbox{Wron}[w_1,v](-S,\lambda)=v(-S,\lambda).
$$
Thus since $w_1'(S,\lambda)$ is known by Lemma \ref{4.1new}, $v(-S,\lambda)$ is also known. Hence 
$$\alpha_n=\frac{\Delta'(\lambda_n)}{v(-S,\lambda_n)},$$
i.e. the norming constants for the Neumann-Neumann problem on $[-S,S]$, given in (\ref{ncfin}), can be uniquely determined.  \qed

To conclude we combine the above theorem together with the results from Section 2 as well as those given in \cite[Sections 3 and 4]{cenw2} to obtain the following Corollary.

\begin{cor}
Given the scattering data i.e. the reflection coefficient and eigenvalues of (\ref{diffapp}), (\ref{tx-conditionapp}) it is possible to obtain the spectra of the Neumann-Neumann and Neumann-Dirichlet problems on $[-S,S]$ or the spectrum and norming constants of the Neumann-Neumann problem on $[-S,S]$. Consequently $M$ can be uniquely reconstructed and $q$ is unique on $\R$.
\end{cor}

{\bf Remark} The fact that, in the above Corollary, we require only the reflection coefficient and eigenvalues (and not the norming constants) of (\ref{diffapp}), (\ref{tx-conditionapp}) in order to uniquely reconstruct $q$ does not contradict the classical results since $q$ has compact essential support, so is known on a very large portion of the line. 
\newsection{Appendix}

The function $v$ obeys the following asymptotics: 

For $\beta =0$ and any arbitrary $\alpha $:

If $0<x\leq S$
\begin{equation}\label{v}
v(x,\lambda ) = \frac{-\sin \sqrt{\lambda }(S-x)}{\sqrt{\lambda }} + O\left(\frac{e^{|\Im \sqrt{\lambda }|(S-x)}}{\lambda }\right)
\end{equation}
\begin{equation}\label{vprime}
v'(x,\lambda ) = \cos \sqrt{\lambda }(S-x) + O\left(\frac{e^{|\Im \sqrt{\lambda }|(S-x)}}{\sqrt{\lambda }}\right).
\end{equation}
For $-S\leq  x<0$ and $m_{12} \not= 0$ 
\begin{equation}\label{vxgzasy}
v(x,\lambda ) = -m_{12} \cos \sqrt{\lambda }S \cos \sqrt{\lambda }x + O\left(\frac{e^{|\Im \sqrt{\lambda }(S-x)|}}{\sqrt{\lambda }}\right)
\end{equation}
\begin{equation}\label{vprimexgzasy}
v'(x,\lambda ) = m_{12} \sqrt{\lambda } \cos \sqrt{\lambda }S \sin \sqrt{\lambda }x + O\left(e^{|\Im \sqrt{\lambda }(S-x)|}\right),
\end{equation}  
if $m_{12} =0$ then
\begin{equation}\label{vxgz2asy}
v(x,\lambda ) = -m_{22} \frac{\sin \sqrt{\lambda }S}{\sqrt{\lambda }}\cos \sqrt{\lambda }x + m_{11}\frac{\sin \sqrt{\lambda }x}{\sqrt{\lambda }}\cos \sqrt{\lambda }S + O\left(\frac{e^{|\Im \sqrt{\lambda }(S-x)|}}{\lambda }\right)
\end{equation}
\begin{equation}\label{vprimexgz2asy}
v'(x,\lambda ) = m_{22}\sin \sqrt{\lambda }S \sin \sqrt{\lambda }x + m_{11} \cos \sqrt{\lambda }S \cos \sqrt{\lambda }x + O\left(\frac{e^{|\Im \sqrt{\lambda }(S-x)|}}{\sqrt{\lambda }}\right).
\end{equation}  

For $\beta \ne 0$ and any arbitrary $\alpha $: 

If $0<x\leq S$  
\begin{eqnarray}
v(x,\lambda) &=& -\sin \beta \cos \sqrt{\lambda} (S-x) +O\left(\frac{e^{|\Im\sqrt{\lambda}|(S-x)}}{\sqrt{\lambda}} \right), \\
v'(x,\lambda) &=& -\sqrt{\lambda}\sin \beta \sin \sqrt{\lambda} (S-x) +O\left(e^{|\Im\sqrt{\lambda}|(S-x)} \right). 
\end{eqnarray}

For $-S\leq  x<0$, if $m_{12}\neq 0$
\begin{eqnarray}
v(x,\lambda) &=& \sqrt{\lambda}m_{12}\sin \beta \cos \sqrt{\lambda}x\sin \sqrt{\lambda} S +O\left(e^{|\Im\sqrt{\lambda}|(S-x)} \right), \label{v1}\\
v'(x,\lambda) &=& -\lambda m_{12}\sin \beta \sin \sqrt{\lambda} x \sin \sqrt{\lambda} S +O\left(\sqrt{\lambda}e^{|\Im\sqrt{\lambda}|(S-x)} \right), \label{v2} 
\end{eqnarray}

if $m_{12}= 0$
\begin{eqnarray}
v(x,\lambda) &=& -m_{22}\sin \beta \cos \sqrt{\lambda}x \cos \sqrt{\lambda} S \nonumber \\
 & & -m_{11}\sin \beta \sin \sqrt{\lambda}x \sin \sqrt{\lambda}S + O\left(\frac{e^{|\Im\sqrt{\lambda}|(S-x)}}{\sqrt{\lambda}} \right), \label{v3}\\
v'(x,\lambda) &=& m_{22}\sqrt{\lambda}\sin \beta \sin \sqrt{\lambda}x \cos \sqrt{\lambda} S \nonumber\\
 && -m_{11}\sqrt{\lambda}\sin \beta \cos \sqrt{\lambda}x \sin \sqrt{\lambda}S + O\left(e^{|\Im\sqrt{\lambda}|(S-x)} \right).\label{v4}
\end{eqnarray}

Working as in the Appendix of \cite{cenw2}, setting $\alpha,\beta\in (0,\pi)$ we have
for $m_{12}=0$
\begin{eqnarray*} 
\Delta_{\alpha,\beta}(\lambda)&=& -\sin \alpha \sin \beta \sqrt{\lambda}(m_{11}+m_{22})\frac{\sin 2 \sqrt{\lambda}S}{2}+O\left(e^{|\Im\sqrt{\lambda}|2S} \right),\\
\Delta_{\alpha,\pi} (\lambda)&  =& 
-\sin\alpha(m_{11}+m_{22})\sin \sqrt{\lambda }S \sin \sqrt{\lambda }S  + O\left(\frac{e^{|\Im \sqrt{\lambda }|2S}}{\sqrt{\lambda }}\right),\\
\Delta_{0,\pi} (\lambda) &=& (m_{22}+ m_{11})\frac{\sin \sqrt{\lambda }S}{\sqrt{\lambda }}\cos \sqrt{\lambda }S + O\left(\frac{e^{|\Im \sqrt{\lambda }|2S}}{\lambda }\right),\\
\Delta_{0,\beta}(\lambda)&=&m_{22}\sin \beta \cos ^2 \sqrt{\lambda}S - m_{11}\sin \beta \sin^2 \sqrt{\lambda}S + O\left(\frac{e^{|\Im\sqrt{\lambda}|2S}}{\sqrt{\lambda}} \right),  
\end{eqnarray*}
and for $m_{12}\neq 0$
\begin{eqnarray*} 
\Delta_{\alpha,\beta}(\lambda)&=&
  \lambda \sin \alpha \sin \beta  m_{12}\sin^2 \sqrt{\lambda}S 
  +O\left(\sqrt{\lambda}e^{|\Im\sqrt{\lambda}|2S} \right),\\
\Delta_{\alpha,\pi} (\lambda) &=& 
  -\sin \alpha m_{12} \sqrt{\lambda } \cos \sqrt{\lambda }S \sin \sqrt{\lambda }S 
  + O\left(e^{|\Im \sqrt{\lambda }|2S}\right),\\
\Delta_{0,\pi} (\lambda) &=& 
  m_{12} \cos^2 \sqrt{\lambda }S  
  + O\left(\frac{e^{|\Im \sqrt{\lambda }|2S}}{\sqrt{\lambda }}\right),\\
\Delta_{0,\beta}(\lambda)&=&-m_{12}\sqrt{\lambda}\sin \beta \cos \sqrt{\lambda}S \sin \sqrt{\lambda}S +O\left(e^{|\Im\sqrt{\lambda}|2S} \right).
\end{eqnarray*}


\end{document}